\documentclass[12pt,a4paper]{article}
\usepackage{lineno,hyperref}
\usepackage[latin5]{inputenc}
\usepackage{amsmath}
\usepackage{graphicx}
\usepackage{subfigure}
\usepackage{float}
\usepackage{mathtools}
\usepackage{color}
\usepackage{cite}
\usepackage{mathrsfs}
\usepackage[]{geometry}
\usepackage[singlelinecheck=false]{caption} % tablo baslıgı sola dayalı
\usepackage{graphicx}
\usepackage{cite}
\usepackage{authblk}
\usepackage{float}
\usepackage{color}
\usepackage{graphicx}
\usepackage{epstopdf}
\usepackage{amscd,amsmath,amssymb,amsthm,amsfonts,epsfig,graphics}

\title{A Comperative Numerical Study Based on Cubic Polynomial and Trigonometric B-splines for the Gardner Equation}
\author{Ozlem Ersoy Hepson$^{a}\thanks{ozersoy@ogu.edu.tr}$, Alper Korkmaz$^{b}$, Idiris Dag$^{c}$  \\
$^{a}${\scriptsize Department of Mathematics \& Computer, Eskisehir Osmangazi University, 26480, Eskisehir, Turkey.}\\
$^{b}${\scriptsize Department of Mathematics, Çankırı Karatekin University, 18200, Çankırı, Turkey.} \\
$^{c}${\scriptsize Department of Computer Engineering, Eskisehir Osmangazi University, 26480, Eskisehir, Turkey.}}
\begin{document}
\maketitle
\begin{abstract} % abstract
Two cubic B-spline functions from different families are placed to the collocation method for the numerical solutions to the Gardner equation. Four models describing propagation of bell shaped single solitary, travel of a kink type wave, wave generation and interaction of two positive bell shaped solitaries propagating in the opposite directions are studied by both methods. The error between the numerical and the analytical solutions is measured by using the discrete maximum norm when the analytical solutions exist. The absolute changes of the lowest three conservation laws are also good indicators of valid results even when the analytical solutions do not exist. The stability of the proposed method is investigated by the Von Neumann analysis.
\end{abstract}
%\begin{keyword}

%\linenumbers
\section{Introduction}
\noindent
The Gardner equation 
\begin{equation}
u_{t}+\alpha uu_{x}+\beta u^2u_{x}+\mu u_{xxx}=0 \label{gardner}
\end{equation}
where $\alpha$, $\beta$ and $\mu$ are nonzero constants is a useful model to study some weakly nonlinear dispersive waves. Having two nonlinear quadratic and cubic terms attracts many researchers to investigate the solutions particulary the numerical ones. Owing to this characteristic of the equation, it is widely known as the KdV-mKdV equation in the literature. The Gardner equation can be completely integrated using the properties of Lax pair and inverse scattering transform. The Miura transformation\cite{miura1} reduces the it to the KdV equation\cite{demler1}. The Gardner equation describes the large amplitude internal waves, too\cite{kamc1}. The implementatiton fields of the equation can be extended to many nonlinear phenomena covering ocean waves, plasma and fluid physics, quantum field theory\cite{wazwazbook}. Using a simple transformation, the Gardner equation can also be converted to a particular form of the modified KdV equation\cite{wazwazbook}. Thus, some solutions including multiple solitons of the equation can be derived from the modified KdV equation\cite{wazwazbook}. The Gardner equation governing the propagation of nonlinear ion-acoustic waves in negative ion loaded plasmas is derived by starting the system of the one dimensional plasma motion equations\cite{ruderman1}. Occurence of rouge waves can also be analyzed by the Gardner equation with the same signed cubic nonlinear and the linear dispersive terms due to allowance modulational instability\cite{grimshaw1}. Diverse family of transcritical flows including kinks, rarefaction waves, classical, reversed and trigonometric undular bores of a stratified fluid over topography are examined in details in the framework of the forced Gardner equation by Kamchatnov et al.\cite{kamc2}.

\noindent
Some particular solutions to model various natural phenomena like interaction of large amplitude solitons or generation are investigated for a particular form having extra first order derivative linear term in the Gardner equation (\ref{gardner})\cite{sul1}. The interactions of waves with different characteristics such as soliton-cnoidal wave, or soliton-periodic wave are derived by the consistent $\tanh$ method\cite{hu1}. The derivation of the soliton-cnoidal wave interaction is based on the consistent Riccati expansion\cite{feng1}. A perturbative study deals with the adiabatic parameter soliton dynamics for the Gardner equation with the aid of perturbation theory\cite{biswas1}. The perturbation parameters for the slow time-dependence of some kink type solitons and correction of order one are developed by using Green's functions\cite{yan1}.

\noindent
So far, various solutions including soliton types, kink-anti kink types, periodic and $N-$soliton solutions to the Gardner equation have been suggested in the related literature. Bekir\cite{bekir1} derives some solutions constructed with $\tanh$ and $\coth$ functions by using the extended $\tanh$ method in which traveling wave transformation and the balance between higher order derivative terms and higher degree power term are the basic keys to the procedure to set up the solutions. In the light of the projective Riccati equations, Fu et al. \cite{fu1} construct many solutions, from solitary waves in terms of hyperbolic functions to periodic ones in terms of sine and cosine functions with the aid of an intermediate transformation in some classical extension method. Combined hyperbolic ansatzes\cite{wazwaz1}, exp-function\cite{akbar1}, $G/G'$-expansion\cite{taghizade1}, improved $G/G'$-expansion\cite{naher1}, generalized $G/G'$-expansion\cite{lu1} and the two-variable $(G'/G,1/G)$ expansion\cite{zayed1} methods are also useful tools to obtain the real and complex exact solutions to the Gardner equation. The dynamics of the solitary wave type solutions of the Gardner equation is examined in details by employing the mapping method\cite{krishan1}. This study also deals with the integration of the perturbed form of the equation using semi inverse principle and gives some solitary waves solutions set up by a particular ansatze method. Jawad\cite{jawad1} also derives an exact solution to the Gardner equation in terms of the tangent function by implementing the $\tan$-$\cot$ function algorithm. Some traveling wave solutions in terms of Weierstrass and Jacobi elliptic functions are constructed in the study of Vassilev et. al.\cite{vasilev1}. Lie group method is applied to develop some new analytical solutions to the Gardner equation\cite{guo1}. A particular similarity reduction is also derived by the Clarkson and Krustal direct method in the same study.  

\noindent
Even though lots of exact and analytical solutions, as summarized in the previous paragraph, are proposed for the Gardner equation, few numerical solutions seem in the literature. A conservative finite difference scheme with the help of the discrete variational methods is implemented to some initial boundary value problems, covering motion of a single soliton and interaction of solitons, for the Gardner equation\cite{nishiyama1}. Conditionally stable restrictive Taylor approximation is another interesting method 
used to simulate kink type solutions for the Gardner equation\cite{rageh1}.

\noindent
This study aims to fill the gap in the related literature by solving some initial boundary value problems for the Gardner equation. We construct polynomial\cite{alp1,alp2} and trigonometric \cite{ctrig} cubic B-spline collocation methods to derive the numerical solutions. Having no parameter inside of both functions makes these two sets different from the exponential cubic B-splines\cite{alp3,oz1,ozsivas,ozkdv}. Various problems covering motion of waves in different profiles, wave generation and interaction of solitary waves are studied numerically. 

\noindent
Having no continuous third order derivatives of the cubic polynomial and trigonometric B-spline functions forces us to reduce the order of the third ordered term in the Gardner Equation. Thus, assuming $v=u_x$ leads the coupled equation system partial differential equations
\begin{equation}
\begin{aligned}
u_{t}+\alpha uu_x+\beta u^{2}u_x+\mu v_{xx}&=0 \\ 
v-u_{x}&=0
\end{aligned}
\label{2}
\end{equation}
The initial data 
\begin{equation}
\begin{aligned}
u(x,0)&=f(x) \\
v(x,0)&=f_x(x)
\end{aligned}
\end{equation}
and the homogenous Neumann boundary conditions are selected from the set
\begin{equation}
\begin{aligned}
u_x(a,t)&=0,\, u_x(b,t)&=0, \\
u_{xx}(a,t)&=0,\, u_{xx}(b,t)&=0, \\
v_x(a,t)&=0,\, v_x(b,t)&=0, \\
v_{xx}(a,t)&=0,\, v_{xx}(b,t)&=0 
\end{aligned}
\end{equation}
in the artifical finite problem interval $[a,b]$.

\section{Methodology}
Consider a uniform equal partition of the problem domain $[x_{0}=a,b=x_{N}]$ with the grids $x_{i}$,$i=0,1,\ldots ,N$ and $h=(b-a)/N$ and the ghost grids $x_{-2}$, $x_{-1}$, $x_{N+1}$, $x_{N+2}$ positioned outside the problem interval. The definition of the polynomial or trigonometric cubic B-splines requires the support of these ghost grids. 

\noindent
The trigonometric cubic B-splines spanning the interval $[a,b]$ are defined as
\begin{equation}
T_{i}(x)=\frac{1}{Z_h }\left \{ 
\begin{tabular}{ll}
$Z^{3}(x_{i-2}),$ & $x\in \left[ x_{i-2},x_{i-1}\right] $ \\ 
$%
Z(x_{i-2})(Z(x_{i-2})\hat{Z}(x_{i})+\hat{Z}(x_{i+1})Z(x_{i-1}))+\hat{Z}(x_{i+2})Z^{2}(x_{i-1}), 
$ & $x\in \left[ x_{i-1},x_{i}\right] $ \\ 
$%
Z(x_{i-2})\hat{Z}^{2}(x_{i+1})+\hat{Z}(x_{i+2})(\hat{Z}(x_{i-1})\hat{Z}(x_{i+1})+\hat{Z}(x_{i+2})Z(x_{i})), 
$ & $x\in \left[ x_{i},x_{i+1}\right] $ \\ 
$\hat{Z}^{3}(x_{i+2}),$ & $x\in \left[ x_{i+1},x_{i+2}\right] $ \\ 
$0,$ & $\text{otherwise}$%
\end{tabular}%
\right.  \label{5a}
\end{equation}%
where $Z(x_{i})=\sin (\frac{x-x_{i}}{2})$, $\hat{Z}(x_{i})=\sin (\frac{x_{i}-x}{2})$, $Z_h =\sin (\frac{h}{2})\sin (h)\sin (\frac{3h}{2})$. The twice continuously differentiable trigonometric cubic B-spline function set $\{T_{i}(x)\}_{i=-1}^{N+1},$ forms a basis for the
functions defined in the same interval \cite{base1,base2}.

\noindent
Similarly, the polynomial cubic B-splines are described as 
\begin{equation}
P_{i}(x)=\frac{1}{h^3}\left \{ 
\begin{tabular}{ll}
$(x-x_{i-2})^3$, & $x\in \left[ x_{i-2},x_{i-1}\right] $ \\ 
$h^3+3h^2(x-x_{i-1})+3h(x-x_{i1})^2-3(x-x_{m-1})^3 $, & $x\in \left[ x_{i-1},x_{i}\right] $ \\ 
$ h^3+3h^2(x_{i+1}-x)+3h(x_{i+1}-x)^2-3(x_{i+1}-x)^3 $, & $x\in \left[ x_{i},x_{i+1}\right] $ \\ 
$(x_{i+2}-x)^3$, & $x\in \left[ x_{i+1},x_{i+2}\right] $ \\ 
$0,$ & $\text{otherwise}$\\
\end{tabular}%
\right.  \label{5b} \\
\end{equation}%
and the set $\{P_{i}(x)\}_{i=-1}^{N+1}$ also constitutes a basis for the functions defined in $[a,b]$.

\noindent
Let $U(x,t)$ and $V(x,t)$ be approximate solutions to $u(x,t)$ and $v(x,t)$ defined as
\begin{equation}
\begin{aligned}
U(x,t)&=\sum_{i=-1}^{N+1}\delta _{i}C_{i}(x)\\
V(x,t)&=\sum_{i=-1}^{N+1}\phi _{i}C_{i}(x) 
\end{aligned}
\label{7}
\end{equation}%
where $C_i(x)=T_{i}(x)$ or $C_i(x)=P_{i}(x)$ and $\delta _{i}$ and $\phi _{i}$ are time dependent parameters that are determined from the collocation points $x_{i},i=0,1, \ldots,N$ and the complementary data. The functional and derivative values of $U(x,t)$ and $V(x,t)$ at a grid $x_i$ is described in terms of time dependent parameters $\delta$ and $\phi$ as
\begin{equation}
\begin{tabular}{l}
$U_{i}=\alpha _{1}\delta _{i-1}+\alpha _{2}\delta _{i}+\alpha _{1}\delta
_{i+1}$ \\ 
$U_{i}^{\prime }=\beta _{1}\delta _{i-1}-\beta _{1}\delta _{i+1}$ \\ 
$U_{i}^{\prime \prime }=\gamma _{1}\delta _{i-1}+\gamma _{2}\delta
_{i}+\gamma _{1}\delta _{i+1}$%
\end{tabular}%
\begin{tabular}{l}
$V_{i}=\alpha _{1}\phi _{i-1}+\alpha _{2}\phi _{i}+\alpha _{1}\phi _{i+1}$
\\ 
$V_{i}^{\prime }=\beta _{1}\phi _{i-1}-\beta _{1}\phi _{i+1}$ \\ 
$V_{i}^{\prime \prime }=\gamma _{1}\phi _{i-1}+\gamma _{2}\phi _{i}+\gamma
_{1}\phi _{i+1}$%
\end{tabular}
\label{8}
\end{equation}%
\noindent
When the trigonometric B-splines are chosen as basis, the coefficients of the time dependent parameters in (\ref{8}) take the forms
\begin{equation}
\begin{array}{lll}
\alpha _{1}=\sin ^{2}(\frac{h}{2})\csc (h)\csc (\frac{3h}{2}) & \beta _{1}=-%
\frac{3}{4}\csc (\frac{3h}{2}) & \gamma _{1}=\dfrac{3((1+3\cos (h))\csc ^{2}(%
\frac{h}{2}))}{16(2\cos (\frac{h}{2})+\cos (\frac{3h}{2}))} \\ 
\alpha _{2}=\dfrac{2}{1+2\cos (h)} & \gamma _{2}=-\dfrac{3\cot ^{2}(\frac{h}{%
2})}{2+4\cos (h)} & 
\end{array}
\label{9a}
\end{equation}
\noindent

Similarly, the same coefficients are determined as
\begin{equation}
\begin{array}{lll}
\alpha _{1}=1 & \beta _{1}=3/h & \gamma _{1}=6/h^2 \\ 
\alpha _{2}=4 & \gamma _{2}=-12/h^2& 
\end{array}
\label{9b}
\end{equation}
when the polynomial cubic B-splines are selected as basis.

\noindent
The Crank-Nicolson and the classical forward finite difference discretization converts the system (\ref{2}) to
\begin{equation}
\begin{array}{r}
\dfrac{U^{n+1}-U^{n}}{\Delta t}+\alpha \dfrac{(UU_x)^{n+1}+(UU_x)^{n}}{2}+\beta\dfrac{(U^{2}U_x)^{n+1}+(U^{2}U_x)^{n}}{2}+\mu\dfrac{%
V_{xx}^{n+1}+V_{xx}^{n}}{2}=0 \\ 
\\ 
\dfrac{V^{n+1}+V^{n}}{2}-\dfrac{U_{x}^{n+1}+U_{x}^{n}}{2}=0%
\end{array}
\label{10}
\end{equation}%
where the superscript $n$ and $n+1$ denotes the functional or derivative values at $n.$th and $n+1.$th time levels, respectively. One should note that $t^{n+1}$ equals $t^{n}+\Delta t$ and $\Delta t$ is the time step size. Substituting the approximate solutions in to the system (\ref{10}) and linearizing the nonlinear terms by the Rubin and Gravis' technique\cite{rubin} as
\begin{equation*}
(UU_x)^{n+1}=U^{n+1}U_x^{n}+U^{n}U_x^{n+1}-U^{n}U_x^{n}
\end{equation*}%
and%
\begin{equation*}
(U^{2}U_x)^{n+1}=2U^{n+1}U^{n}U_x^{n}+(U^{n})^{2}U_x^{n+1}-2(U^{n})^{2}U_x^{n}
\end{equation*}
yields
\scriptsize
\begin{equation}
\begin{aligned}
&\left[ \left( \frac{2}{\Delta t}+\alpha L+2\beta KL\right) \alpha _{1}+\left( \alpha K+\beta K^{2}\right)
\beta _{1}
\right] \delta _{j-1}^{n+1}+\left[ \mu\gamma _{1}\right] \phi _{j-1}^{n+1}+\left[ \left( \frac{%
2}{\Delta t}+\alpha L+2\beta KL\right) \alpha _{2}\right] \delta
_{j}^{n+1}\\
&+\left[ \mu_{3}\gamma _{2}\right] \phi _{j}^{n+1} \left[ \left( \frac{2}{\Delta t}+\alpha L+2\beta KL\right) \alpha _{1}-\left( \alpha K+\beta K^{2}\right)
\beta _{1}\right] \delta _{j+1}^{n+1}+\left[ \mu\gamma _{1}\right] \phi _{j+1}^{n+1}\\
&=\left[ \left( \frac{2}{\Delta t}+\beta KL\right) \alpha _{1}\right]
\delta _{j-1}^{n}-\mu\gamma _{1}\phi _{j-1}^{n}+\left[ \left( \frac{2}{%
\Delta t}+\beta KL\right) \alpha _{2}\right] \delta _{j}^{n}-\mu
_{3}\gamma _{2}\phi _{j}^{n}+\left[ \left( \frac{2}{\Delta t}+\mu
_{2}KL\right) \alpha _{1}\right] \delta _{j+1}^{n}&-\mu\gamma _{1}\phi
_{j+1}^{n} \\
&-\beta _{1}\delta _{j-1}^{n+1}+\alpha _{1}\phi _{j-1}^{n+1}+\alpha _{2}\phi
_{j}^{n+1}+\beta _{1}\delta _{j+1}^{n+1}+\alpha _{1}\phi _{j+1}^{n+1}
=\beta _{1}\delta _{j-1}^{n}-\alpha _{1}\phi _{j-1}^{n}-\alpha _{2}\phi
_{j}^{n}-\beta _{1}\delta _{j+1}^{n}-\alpha _{1}\phi _{j+1}^{n} \\
&m=0,...,N,\quad n=0,1...,
\end{aligned}\label{12}
\end{equation}

\normalsize
where%
\begin{equation*}
\begin{aligned}
L_{1}&=\alpha _{1}\delta _{i-1}^{n}+\alpha _{2}\delta _{i}^{n}+\alpha _{1}\delta
_{i+1}^{n} \\ 
L_{2}&=\alpha _{1}\phi _{i-1}^{n}+\alpha _{2}\phi _{i}^{n}+\alpha _{1}\phi
_{i+1}^{n}%
\end{aligned}%
\end{equation*}%

\noindent
For the simplicity, we use the matrix notation
\begin{equation}
\mathbf{Ad}^{n+1}=\mathbf{Bd}^{n}  \label{13}
\end{equation}%

with the matrix elements
\begin{equation*}
\mathbf{A=}%
\begin{matrix}
\eta _{1} & \eta _{2} & \eta _{3} & \eta _{4} & \eta _{5} & \eta _{2} &  & 
&  &  \\ 
-\beta _{1} & \alpha _{1} & 0 & \alpha _{2} & \beta _{1} & \alpha _{1} &  & 
&  &  \\ 
&  & \eta _{1} & \eta _{2} & \eta _{3} & \eta _{4} & \eta _{5} & \eta _{2} & 
&  \\ 
&  & -\beta _{1} & \alpha _{1} & 0 & \alpha _{2} & \beta _{1} & \alpha _{1}
&  &  \\ 
&  &  & \ddots  & \ddots  & \ddots  & \ddots  & \ddots  & \ddots  &  \\ 
&  &  &  & \eta _{1} & \eta _{2} & \eta _{3} & \eta _{4} & \eta _{5} & \eta
_{2} \\ 
&  &  &  & -\beta _{1} & \alpha _{1} & 0 & \alpha _{2} & \beta _{1} & \alpha
_{1}%
\end{matrix}%
\end{equation*}

\begin{equation*}
\mathbf{B=}%
\begin{matrix}
\eta _{6} & -\eta _{2} & \eta _{7} & -\eta _{4} & \eta _{6} & -\eta _{2} &  & 
&  &  \\ 
\beta _{1} & -\alpha _{1} & 0 & -\alpha _{2} & -\beta _{1} & -\alpha _{1} & 
&  &  &  \\ 
&  & \eta _{6} & -\eta _{2} & \eta _{7} & -\eta _{4} & \eta _{6} & -\eta _{2} & 
&  \\ 
&  & \beta _{1} & -\alpha _{1} & 0 & -\alpha _{2} & -\beta _{1} & -\alpha
_{1} &  &  \\ 
&  &  & \ddots & \ddots & \ddots & \ddots & \ddots & \ddots &  \\ 
&  &  &  & \eta _{6} & -\eta _{2} & \eta _{7} & -\eta _{4} & \eta _{6} & -\eta _{2} \\ 
&  &  &  & \beta _{1} & -\alpha _{1} & 0 & -\alpha _{2} & -\beta _{1} & 
-\alpha _{1}%
\end{matrix}%
\end{equation*}%
and%
\begin{equation*}
\begin{array}{ll}
\eta _{1}=\left( \frac{2}{\Delta t}+\mu _{1}L+2\mu _{2}KL\right) \alpha _{1}+\left( \mu _{1}K+\mu _{2}K^{2}\right)
\beta _{1}
& \eta _{6}=\left( \frac{2}{\Delta t}+\mu _{2}KL\right) \alpha _{1} \\ 
\eta _{2}=\mu _{3}\gamma_{1} &  \eta _{5}=\left( \frac{2}{\Delta t}+\mu _{1}L+2\mu _{2}KL\right) \alpha _{2}-\left( \mu _{1}K+\mu _{2}K^{2}\right)
\beta _{1}\\ 
\eta _{3}=\left( \frac{2}{\Delta t}+\mu _{1}L+2\mu _{2}KL\right) \alpha _{2}-\left( \mu _{1}K+\mu _{2}K^{2}\right)
\beta _{1}
& \eta _{7}=\left( \frac{2}{\Delta t}+\mu _{2}KL\right) \alpha _{2} \\ 
\eta _{4}=\mu _{3}\gamma
_{2} & 
\end{array}%
\end{equation*}

Even though there are $2N+6$ unknown
parameters 
\begin{equation*}
\mathbf{d}^{n+1}=(\delta _{-1}^{n+1},\phi _{-1}^{n+1},\delta _{0}^{n+1},\phi
_{0}^{n+1},\ldots ,\delta _{n+1}^{n+1},\phi _{n+1}^{n+1},).
\end{equation*}
in the system (\ref{13}), it consists of only $2N+2$ linear equations. In order to equalize the numbers of equations and unknowns, we eliminate some unknowns by manipulating the appropriate boundary conditions. 

%The homogenous Neumann conditions $U_{x}(a,t)=0$, $V_{x}(a,t)=0$ and $U_{x}(b,t)=0$, $V_{x}(b,t)=0$ give the relation between the unknown parameters as
%\begin{equation*}
%\begin{array}{l}
%\delta _{-1}^{n+1}=\delta _{1}^{n+1} \\ 
%\phi _{-1}^{n+1}=\phi _{1}^{n+1} \\ 
%\delta _{N-1}^{n+1}=\delta _{N+1}^{n+1} \\ 
%\phi _{N-1}^{n+1}=\phi _{N+1}^{n+1}
%\end{array}
%\end{equation*}

%{\color{red}\textbf{
%İkinci Türev koşullarının eliminasyonda kullanılması nasıl olacak????????????????}
%}

\noindent
%Thus, eliminating the parameters $\delta _{-1}^{n+1},\phi _{-1}^{n+1},\delta _{N+1}^{n+1},\phi
%_{N+1}^{n+1}$ from the system (\ref{13}), gives a solvable system having $2N+2$ unknowns and equations. In order to start the iteration, we should determine the initial vector $\mathbf{x^0}$. Some algebraic manipulations of initial and boundary data give 
%\begin{equation*}
%\begin{array}{l}
%U_{x}(a,0)=0=\delta _{-1}^{0}-\delta _{1}^{0}, \\ 
%U(x_{i},0)=\alpha _{1}\delta _{i-1}^{0}+\alpha _{2}\delta _{i}^{0}+\alpha
%_{1}\delta _{i+1}^{0}=u(x_{i},0),i=1,...,N-1 \\ 
%U_{x}(b,0)=0=\delta _{N-1}^{0}-\delta _{N+1}^{0}, \\ 
%V_{x}(a,0)=0=\phi _{-1}^{0}-\phi _{1}^{0} \\ 
%V(x_{i},0)=\alpha _{1}\phi _{i-1}^{0}+\alpha _{2}\phi _{i}^{0}+\alpha
%_{1}\phi _{i+1}^{0}=v(x_{i},0),i=1,...,N-1 \\ 
%V_{x}(a,0)=\phi _{N-1}^{0}-\phi _{N+1}^{0}%
%\end{array}%
%\end{equation*}
\noindent
%Since the initial vector can be calculated by using the equations given above, the iteration algorithm is ready to run.

\section{Stability Analysis}
The stability of the method is investigated by performing the Von-Neumann analysis where 
\begin{eqnarray}
\delta_{j}^{n} &=&K_1\xi ^{n}\exp (ij\varphi )  \label{k} \\
\phi _{j}^{n} &=&K_2\xi ^{n}\exp (ij\varphi )  \notag
\end{eqnarray}%
\begin{equation*}
\rho=\frac{\xi ^{n+1}}{\xi ^{n}}
\end{equation*}%
Here, $K_1$ and $K_2$ represent the harmonics amplitude, $k$ is the
mode number, $\rho$ is the amplification factor and $\varphi =kh$. The term $U+U^2$ is assumed as locally constant and replaced by $\varepsilon$. Substituting \ref{k} into the system 
\begin{eqnarray}
&&\alpha_{1}\delta_{j-1}^{n+1}+\alpha_{2}\delta_{j}^{n+1}+\alpha_{1}\delta_{j+1}^{n+1}+\frac{\lambda
k\varepsilon }{2}(\beta_{1}\delta_{j-1}^{n+1}-\beta_{1}\delta_{j+1}^{n+1})+\frac{k\mu }{2}%
(\gamma_{1}\phi _{j-1}^{n+1}+\gamma_{2}\phi _{j}^{n+1}+\gamma_{1}\phi _{j+1}^{n+1})
\label{k1} \\
&=&\alpha_{1}\delta_{j-1}^{n}+\alpha_{2}\delta_{j}^{n}+\alpha_{1}\delta_{j+1}^{n}-\frac{\lambda
k\varepsilon }{2}(\beta_{1}\delta_{j-1}^{n}-\beta_{1}\delta_{j+1}^{n})-\frac{k\mu }{2}%
(\gamma_{1}\phi _{j-1}^{n}+\gamma_{2}\phi _{j}^{n}+\gamma_{1}\phi _{j+1}^{n})  \notag
\end{eqnarray}%
\begin{eqnarray}
&&\beta_{1}\delta_{j-1}^{n+1}-\beta_{1}\delta_{j+1}^{n+1}-\alpha_{1}\phi _{j-1}^{n+1}-\alpha_{2}\phi
_{j}^{n+1}-\alpha_{1}\phi _{j+1}^{n+1}  \label{k2} \\
&=&-\beta_{1}\delta_{j-1}^{n}+\beta_{1}\delta_{j+1}^{n}+\alpha_{1}\phi _{j-1}^{n}+\alpha_{2}\phi
_{j}^{n}+\alpha_{1}\phi _{j+1}^{n}  \notag
\end{eqnarray}%
gives
\begin{eqnarray*}
&&\xi ^{n+1}\left[ K_1\left( 2\alpha_{1}\cos \varphi +\alpha_{2}\right) +\frac{K_2k\mu }{2}%
\left( 2\gamma_{1}\cos \varphi +\gamma_{2}\right) -i\lambda k\varepsilon K_1\beta_{1}\sin
\varphi \right]  \\
&=&\xi ^{n}\left[ K_1\left( 2\alpha_{1}\cos \varphi +\alpha_{2}\right) -\frac{K_2k\mu }{2}%
\left( 2\gamma_{1}\cos \varphi +\gamma_{2}\right) +i\lambda k\varepsilon K_1\beta_{1}\sin
\varphi \right] 
\end{eqnarray*}%
\begin{equation*}
\frac{\xi ^{n+1}}{\xi ^{n}}=\frac{\left[ K_1\left( 2\alpha_{1}\cos \varphi
+\alpha_{2}\right) -\frac{K_2k\mu }{2}\left( 2\gamma_{1}\cos \varphi +\gamma_{2}\right)
+i\lambda k\varepsilon K_1\beta_{1}\sin \varphi \right] }{\left[ K_1\left(
2\alpha_{1}\cos \varphi +\alpha_{2}\right) +\frac{K_2k\mu }{2}\left( 2\gamma_{1}\cos \varphi
+\gamma_{2}\right) -i\lambda k\varepsilon K_1\beta_{1}\sin \varphi \right] }
\end{equation*}%
\begin{equation}
\rho=\frac{\xi ^{n+1}}{\xi ^{n}}=\frac{M_{1}+iK_{1}}{M_{2}-iK_{1}}  \label{k3}
\end{equation}%
where%
\begin{eqnarray*}
M_{1} &=&K_1\left( 2\alpha_{1}\cos \varphi +\alpha_{2}\right) -\frac{K_2k\mu }{2}\left(
2\gamma_{1}\cos \varphi +\gamma_{2}\right)  \\
M_{2} &=&K_1\left( 2\alpha_{1}\cos \varphi +\alpha_{2}\right) +\frac{K_2k\mu }{2}\left(
2\gamma_{1}\cos \varphi +\gamma_{2}\right)  \\
K_{1} &=&\lambda k\varepsilon K_1\beta_{1}\sin \varphi 
\end{eqnarray*}%
and
\begin{eqnarray*}
&&\xi ^{n+1}\left[ -K_2\left( 2\alpha_{1}\cos \varphi +\alpha_{2}\right) -2iK_1\beta_{1}\sin
\varphi \right]  \\
&=&\xi ^{n}\left[ K_2\left( 2\alpha_{1}\cos \varphi +\alpha_{2}\right) +2iK_1\beta_{1}\sin
\varphi \right] 
\end{eqnarray*}%
\begin{equation*}
\frac{\xi ^{n+1}}{\xi ^{n}}=\frac{K_2\left( 2\alpha_{1}\cos \varphi +\alpha_{2}\right)
+2iK_1\beta_{1}\sin \varphi }{-K_2\left( 2\alpha_{1}\cos \varphi +\alpha_{2}\right)
-2iK_1\beta_{1}\sin \varphi }
\end{equation*}%
\begin{equation}
\rho =\frac{\xi ^{n+1}}{\xi ^{n}}=\frac{M_{3}+iK_{2}}{M_{4}-iK_{2}}  \label{k4}
\end{equation}%
\begin{eqnarray*}
M_{3} &=&K_2\left( 2\alpha_{1}\cos \varphi +\alpha_{2}\right)  \\
M_{4} &=&-K_2\left( 2\alpha_{1}\cos \varphi +\alpha_{2}\right)  \\
K_{2} &=&2K_1\beta_{1}\sin \varphi 
\end{eqnarray*}
It can be concluded from both (\ref{k3}) and (\ref{k4}) that $\left \vert \rho \right \vert$ is less than or equal to $1$. Thus, the proposed method method is unconditionally stable.

\section{Numerical Examples}
\noindent
In this section, we solve some analytical and non-analytical initial boundary value problems to validate the proposed methods and compare the results produced by both methods. We give some tools below to measure the accuracies of the proposed methods for a fair comparison.

\noindent
The discrete maximum norm measuring the error between the numerical and the analytical solutions, if exist, defined as
\begin{equation*}
\begin{aligned}
L_{\infty}(t)&=\left \vert u(x,t)-U(x,t)\right \vert _{\infty }=\max \limits_{i}\left
\vert u(x_i,t)-U(x_i,t)\right \vert \\
\end{aligned}%
\end{equation*}
at the time $t$ is chosen as a confirmative for the accuracy and validity of the proposed method. The conservation laws can also be good indicators of the validity even when the analytical solutions do not exist. The lowest three conservation laws describing the momentum(M), the energy(E) and the Hamiltonian(H) for the Gardner equation are defined as
\begin{equation}
\begin{aligned}
M&=\int\limits_{-\infty}^{\infty}{udx} \\
E&=\int\limits_{-\infty}^{\infty}{u^2dx} \\
H&=\int\limits_{-\infty}^{\infty}{\dfrac{\alpha u^3}{3}+\dfrac{\beta u^4}{6}-\mu (u_x)^2dx}
\end{aligned}
\end{equation} 
and are expected to remain constant as time goes\cite{hamdi1}. Define the absolute relative changes $C(M_t)$, $C(E_t)$ and $C(H_t)$ of the conservation laws $M$, $E$ and $H$ as
\begin{equation}
\begin{aligned}
C(M_t)&=\left | \frac{M_t-M_0}{M_0} \right | \\
C(E_t)&=\left | \frac{E_t-E_0}{E_0} \right | \\
C(H_t)&=\left| \frac{H_t-H_0}{H_0} \right|
\end{aligned}
\end{equation}
where $M_0$, $E_0$ and $H_0$ are initial, $M_t$, $E_t$ and $H_t$ are the numerically computed values of the conserved quantities at the time $t$.

\subsection{Propagation of Bell Type Solitary Wave}
\begin{figure}[hp]
    \subfigure[Propagation of the bell shape solitary]{
   \includegraphics[scale =0.6] {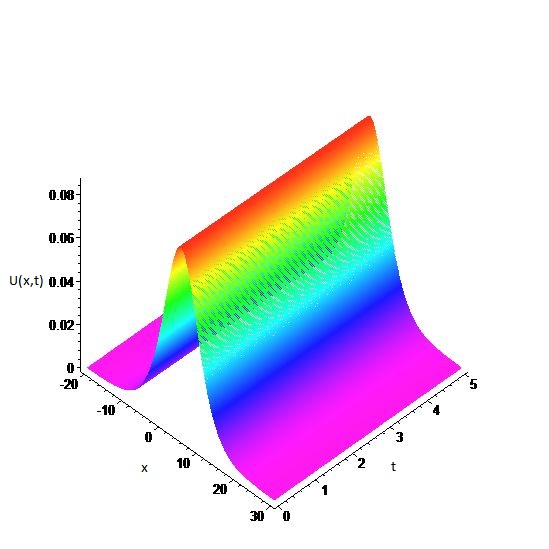}
   \label{fig:1a}
 }
 \subfigure[The maximum error distribution for polynomial cubic B-spline method]{
   \includegraphics[scale =0.6] {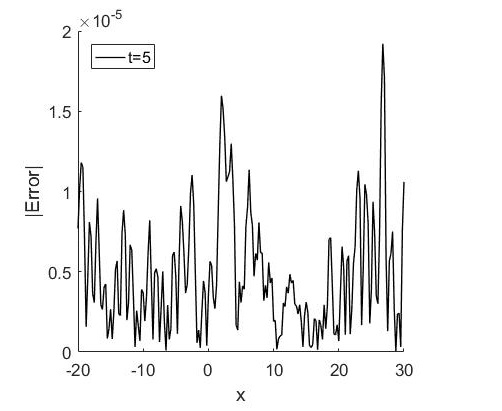}
   \label{fig:1b}
 }
 \subfigure[The maximum error distribution for trigonometric cubic B-spline method]{
   \includegraphics[scale =0.6] {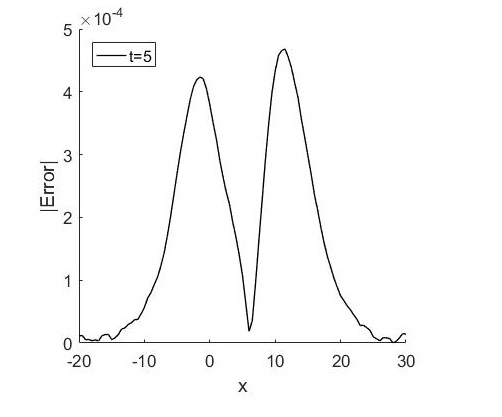}
   \label{fig:1c}
 }
 \caption{Propagation of the bell shape solitary and the maximum error distributions at $t=5$ with $N=100$}
\end{figure}

In the Gardner equation,we assume that the parameters $\alpha=4,$ $\beta=-3$ and $\mu=1$. Thus, the analytical solution of the Gardner equation takes the form
\begin{equation}
u(x,t)=\dfrac{2}{12+3\sqrt{14}\cosh{\left(-\frac{x}{3}+\frac{5}{3}+\frac{t}{27}\right)}} \label{sech}
\end{equation}
\noindent
This solution is adapted from the $\cosh$ ansatze solution describing the propagation of a positive bell type solitary wave along the $x$-axis\cite{wazwaz1}. The initial condition is derived by substituting $t=0$ into the analytical solution. The first order homogeneous Neumann conditions are applied at both end of the finite interval $[-20,30]$. For the sake of comparison, the proposed routines are run up to the terminating time $t=5$ with the fixed discretization parameter $\Delta t=0.1$ and various numbers of grids $N$. A three dimensional simulation of the propagation, Fig \ref{fig:1a} and the maximum error distributions of the solution obtained by using both methods at the simulation terminating time are depicted in Fig \ref{fig:1b}-\ref{fig:1c}. 

\noindent
The accumulation of the error around the peak position draws the attention in the trigonometric method as no clear accumulation interval for the error in the polynomial method. The comparison of the discrete maximum norms at the half and the end times of the simulation process are tabulated in the Table \ref{t1}. The error is in five decimal digits in polynomial method but is in four decimal digits in trigonometric method at the simulation terminating time $t=5$ for $N=100$. The increase of the number of grids to $200$ decreases the value of the norm but the errors are still in the same digits as the previous choice of $N$. Increasing the number of points to $300$ improves the results of the trigonometric method more than the improvement of the polynomial method. Finally, five decimal digit accuracy of the results is obtained for both methods when $N$ is chosen as $400$. 

\begin{table}[hp]
\caption{Comparison of the discrete maximum norms} 
	\begin{tabular}{lllll}
			\hline\hline
			&Polynomial method && Trigonometric method \\
			$N$ & $L_{\infty }(2.5)$ & $L_{\infty }(5)$ & $L_{\infty }(2.5)$ & 
$L_{\infty }(5)$ \\ 	\hline
$100$ & $3.2726\times 10^{-5}$ & $5.2261\times 10^{-5}$ & $2.3233\times
10^{-4}$ & $4.6808\times 10^{-4}$ \\ 
$200$ & $2.0539\times 10^{-5}$ & $1.9207\times 10^{-5}$ & $6.6360\times
10^{-5}$ & $1.1664\times 10^{-4}$ \\ 
$300$ & $1.4433\times 10^{-5}$ & $1.7145\times 10^{-5}$ & $3.5534\times
10^{-5}$ & $5.7732\times 10^{-5}$ \\ 
$400$ & $1.4609\times 10^{-5}$ & $1.6283\times 10^{-5}$ & $2.3168\times
10^{-5}$ & $3.5892\times 10^{-5}$ \\ \hline \hline
	\end{tabular}
	\label{t1}
\end{table}

\noindent
The initial values of the lowest three conservation laws are determined analytically as
\begin{equation}
\begin{aligned}
M_0&= 4\,{\it arctanh} \left(  A^{-1}
 \right) \sqrt {2} \\
E_0&=4/3\,{\frac {16\,{\it arctanh} \left(  A ^{-1} \right) +4\,{\it arctanh} \left(  A ^{-1} \right) \sqrt {2}\sqrt {7}-A}{A}} \\
H_0&= {\frac {4}{81}}\,{\frac {4833024\,{\it arctanh} \left(  A ^{-1} \right) -629300\,\sqrt {2}-336375\,
\sqrt {7}+1291680\,{\it arctanh} \left(  A^{-1} \right) \sqrt {2}\sqrt {7}}{13455\,\sqrt {7}+25172\,
\sqrt {2}}}
\end{aligned}
\end{equation}
where $A=2\,\sqrt {2}+\sqrt {7}$ in the infinite interval $(-\infty,\infty)$. The approximate values of these quantities are calculated as $M_0=1.045100915$, $E_0=0.06013455349$ and $H_0=0.004070220312$ at the beginning of the numerical simulation and are expected to remain constant as time goes. The absolute relative changes of the conservation laws are reported in Table \ref{t2}. Both methods can be said to generate successful results owing to the absolute relative changes reported in the table but it is not easy to decide which one gives more accurate results than the other one.
 
\begin{table}[hp]
\caption{The absolute relative changes of the conservation laws for the bell type solitary} 
	\begin{tabular}{llll}
			\hline\hline	
$N$ &  $C(M_{5})$ & $C(E_{5})$ & $C(H_{5})$  \\ 
			\hline\hline	
		& Polynomial Method &&\\
			\hline
$100$ & $5.4401\times 10^{-6}$ & $3.8347\times 10^{-8}$ & $1.5235\times 10^{-6}$ \\
$200$ &  $3.4007\times 10^{-6}$ & $5.1804\times 10^{-8}$ & $1.6997\times 10^{-6}$ \\
$300$ & $6.1231\times 10^{-8}$ & $2.3335\times 10^{-8}$ & $2.8345\times 10^{-6}$ \\
$400$ & $8.2790\times 10^{-7}$ & $2.6775\times 10^{-9}$ & $3.3932\times 10^{-6}$ \\ 

& Trigonometric Method &&\\
\hline
$100$ & $2.1937\times 10^{-6}$ & $1.4472\times 10^{-8}$ & $7.0199\times 10^{-5}$ \\ 
$200$ & $3.0848\times 10^{-6}$ & $4.6159\times 10^{-8}$ & $2.5761\times 10^{-6}$ \\ 
$300$ &  $8.2228\times 10^{-8}$ & $2.0955\times 10^{-8}$ & $1.9806\times 10^{-6}$ \\ 
$400$ & $1.3187\times 10^{-6}$ & $6.8085\times 10^{-10}$ & $3.1222\times 10^{-6}$ \\ 
\hline
 \hline 
	\end{tabular}
	\label{t2}
\end{table}

\subsection{Motion of Kink Type Wave}
\begin{figure}[hp]
    \subfigure[Motion of the kink type wave ]{
   \includegraphics[scale =0.6] {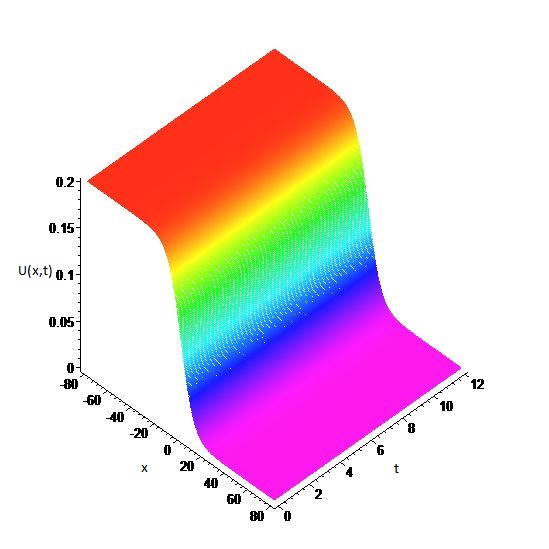}
   \label{fig:2a}
 }
 \subfigure[The maximum error distribution for polynomial cubic B-spline method]{
   \includegraphics[scale =0.6] {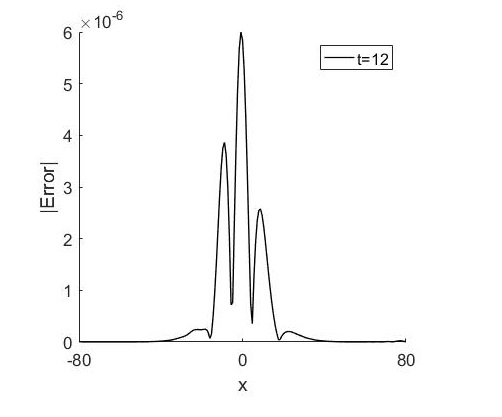}
   \label{fig:2b}
 }
 \subfigure[The maximum error distribution for trigonometric cubic B-spline method]{
   \includegraphics[scale =0.7] {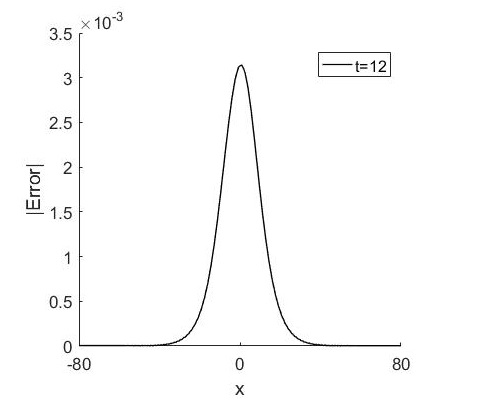}
   \label{fig:2c}
 }
 \caption{Motion of the kink type wave and the maximum error distributions at $t=12$ with $N=200$}
\end{figure}

\noindent
Kink type wave solution of the Gardner equation is 
\begin{equation}
u(x,t)=\dfrac{1}{10}-\dfrac{1}{10}\tanh{(\dfrac{\sqrt{30}}{60}(x-\dfrac{1}{30}t))} \label{tanh}
\end{equation}
for the coefficient choice $\alpha=1$, $\beta=-5$ and $\mu=1$\cite{wazwazbook}. This kink type wave travels to the right with the velocity $1/30$, Fig \ref{fig:2a}. The the initial condition required to start the iteration of the time integration is determined by assuming $t=0$ in the analytical solution (\ref{tanh}). We choose homogeneous second order and first order Neumann conditions in the polynomial and trigonometric methods, respectively. Both routines are run up to the terminating time $t=12$ for various selections of discretization parameter $N$ with a fixed time step size $\Delta t=0.1$ in the finite interval $[-80,80]$. The error distributions for both the polynomial and the trigonometric methods are plotted in Fig \ref{fig:2b}-\ref{fig:2c}. The errors are observed to be accumulated around the descent positions of the wave as expected.

\noindent
A comparison of both methods are examined in details for various numbers of grid points in Table \ref{t3}. When $N=100$, the maximum error is in five decimal digits in the results of the polynomial method but the results of trigonometric method have 2 decimal digits accuracy at the end of the simulation. Increasing $N$ to $200$ improves the accuracy of results of both methods one more decimal digits. $N=400$ gives six decimal digit accuracy in the polynomial method as four decimal digits in the trigonometric method. When the number of grids are chosen as $600$, the results to seven decimal digits in the polynomial method but we do not observe an improvement of the results in decimal digits in the trigonometric method although the results are improved almost two times. The final experiment is achieved by using $N=800$. The results are improved in both methods but not sufficient to change the decimal digits accuracies. 

\begin{table}[hp]
\caption{Comparison of the discrete maximum norms for the kink wave motion} 
	\begin{tabular}{llll}
			\hline\hline
			$N$ & $L_{\infty }(4)$ & $L_{\infty }(8)$ & $L_{\infty }(12)$\\ 
			\hline
			&Polynomial method  \\
\hline
$100$ & $8.4150\times 10^{-6}$ & $1.6152\times 10^{-5}$ & $2.3157\times 10^{-5}$  \\ 
$200$ & $2.1204\times 10^{-6}$ & $4.1172\times 10^{-6}$ & $5.9958\times 10^{-6}$ \\ 
$400$ & $5.3320\times 10^{-7}$ & $1.0346\times 10^{-6}$ & $1.5016\times 10^{-6}$ \\ 
$600$ & $2.3698\times 10^{-7}$ & $4.6025\times 10^{-7}$ & $6.6857\times 10^{-7}$ \\ 
$800$ & $1.3337\times 10^{-7}$ & $2.5929\times 10^{-7}$ & $3.7317\times 10^{-7}$  \\
\hline 
&Trigonometric method \\
\hline
$100$ &  $4.4735\times 10^{-3}$ & $8.9512\times 10^{-3}$ & $1.3436\times 10^{-2}$ \\ 
$200$  & $1.0441\times 10^{-3}$ & $2.0912\times 10^{-3}$ & $3.1431\times 10^{-3}$ \\ 
$400$  & $2.6140\times 10^{-4}$ & $5.2439\times 10^{-4}$ & $7.8804\times 10^{-4}$ \\ 
$600$  & $1.2642\times 10^{-4}$ & $2.5465\times 10^{-4}$ & $3.8356\times 10^{-4}$ \\ 
$800$  & $8.7432\times 10^{-5}$ & $1.7743\times 10^{-4}$ & $2.6841\times 10^{-4}$ \\
 \hline \hline
	\end{tabular}
	\label{t3}
\end{table}

\noindent
The initial analytical values of the conservation laws are computed in the problem interval $[-80,80]$ via symbolic computation software as
\begin{equation}
{\scriptsize
\begin{aligned}
M_0&=16 \\
E_0&= -1/25\,{\frac {-\sqrt {30}-80\,{{\rm e}^{8/3\,\sqrt {30}}}+\sqrt {30}{
{\rm e}^{8/3\,\sqrt {30}}}-80}{1+{{\rm e}^{8/3\,\sqrt {30}}}}}\cong 2.980911178 \\
H_0&= {\frac {1}{1125}}\,{\frac { \left( -2\,\sqrt {30}+120 \right) {{\rm e}
^{8\,\sqrt {30}}}+ \left( 360+3\,\sqrt {30} \right) {{\rm e}^{8/3\,
\sqrt {30}}}+ \left( -3\,\sqrt {30}+360 \right) {{\rm e}^{16/3\,\sqrt 
{30}}}+2\,\sqrt {30}+120}{1+3\,{{\rm e}^{8/3\,\sqrt {30}}}+3\,{{\rm e}
^{16/3\,\sqrt {30}}}+{{\rm e}^{8\,\sqrt {30}}}}} \cong 0.09692938338
\end{aligned}
}
\end{equation} 
\noindent
The conservation laws are recorded at the end of the simulation time $t=12$ for all choices of number of grids, Table \ref{t4}. The absolute relative changes of conservation laws seem insensitive to the increase the number of the grids for the polynomial method. On the other hand, increasing the number of grids from $100$ to $200$ improves the sensitivity in decimal digits for the trigonometric method. To keep to increase the number of grids decreases the decimal digit sensitivity of the conservation laws.

\begin{table}[hp]
\caption{The absolute relative changes of the conservation laws for the travel of the kink type wave} 
	\begin{tabular}{llll}
			\hline\hline	
$N$ &  $C(M_{12})$ & $C(E_{12})$ & $C(H_{12})$  \\ 
			\hline\hline	
		& Polynomial Method &&\\
			\hline
$100$  & $4.9504 \times 10^{-3}$ &$ 
5.3104\times 10^{-3}$ &$ 5.4423\times 10^{-3}$  \\ 
$200$  & $4.9751\times 10^{-3}$ &$ 
5.3388\times 10^{-3} $&$ 5.4721\times 10^{-3}$ \\ 
$400$ & $4.9875\times 10^{-3}$ & $
5.3531\times 10^{-3}$ & $5.4871\times 10^{-3}$  \\ 
$600$  & $4.9916\times 10^{-3}$ & $
5.3578\times 10^{-3} $& $5.4922\times 10^{-3} $\\ 
$800$  & $4.9937\times 10^{-3} $& $
5.3602\times 10^{-3} $&$ 5.4947\times 10^{-3}$ \\ \hline
& Trigonometric Method &&\\
\hline
$100$  &$ 1.4566\times 10^{-2} $&$ 
1.5625\times 10^{-2} $&$ 1.6014\times 10^{-2}$ \\ 
$200$& $3.7506\times 10^{-4}$ &$ 
4.0270\times 10^{-4}$ & $4.1256\times 10^{-4}$ \\ 
$400$  & $3.8575\times 10^{-3}$ & $
4.1408\times 10^{-3} $& $4.2444\times 10^{-3}$ \\ 
$600$  & $4.4909\times 10^{-3} $&$ 
4.8213\times 10^{-3} $&$ 4.9421\times 10^{-3}$ \\ 
$800$  &$ 4.7123\times 10^{-3}$ &$ 
5.0596\times 10^{-3} $& $5.1862\times 10^{-3}$ \\ 
\hline\hline
	\end{tabular}
	\label{t4}
\end{table}

\subsection{Wave Generation from an Initial Pulse}
\noindent
The perturbed Gardner equation of the form
\begin{equation}
u_{t}+\alpha uu_{x}+\beta u^2u_{x}+\mu u_{xxx}=\epsilon \label{gardnerp}
\end{equation}
for some nonzero real $\epsilon$ can be useful to study the wave generation from an initial positive pulse. The decomposition of the balance among the nonlinear terms and the third order derivative is expected not to keep the shape or velocity as propagating. Thus, the initial condition is generated from the initial condition of the first problem sensitively as
\begin{equation}
u(x,t)=\dfrac{2}{3}\dfrac{1}{4+\sqrt{14}\cosh{\left(\dfrac{x}{3}-\dfrac{5}{3}\right)}}
\end{equation}
by perturbating the initial condition. We assume that $\alpha=10$, $\beta=-3$ and $\mu=1$ in the Gardner equation (\ref{gardner}). We run the proposed algorithms with the discretization parameters $N=400$ and $\Delta t=0.01$ in the artificial problem interval $[-40,60]$ up to the time $t=15$. 

\begin{figure}[hp]
    \subfigure[Initial data]{
   \includegraphics[scale =0.65] {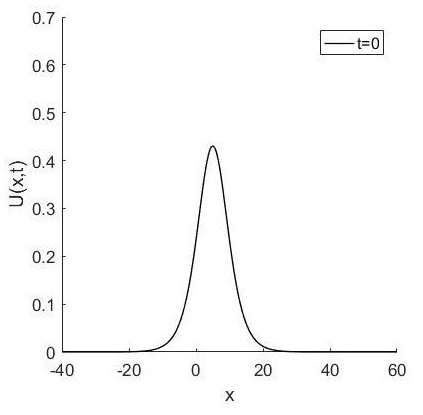}
   \label{fig:3a}
 }
   \subfigure[$t=5$]{
   \includegraphics[scale =0.65] {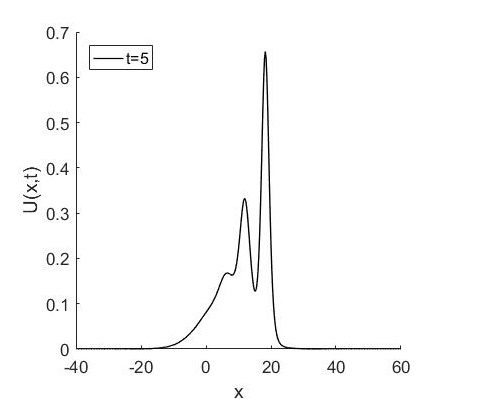}
   \label{fig:3b}
 }
 \subfigure[$t=10$]{
   \includegraphics[scale =0.65] {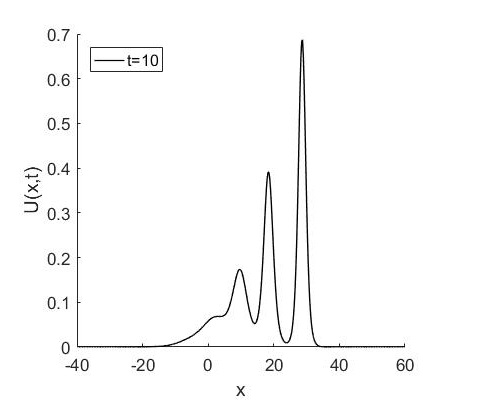}
   \label{fig:3c}
 }
 \subfigure[$t=15$]{
   \includegraphics[scale =0.65] {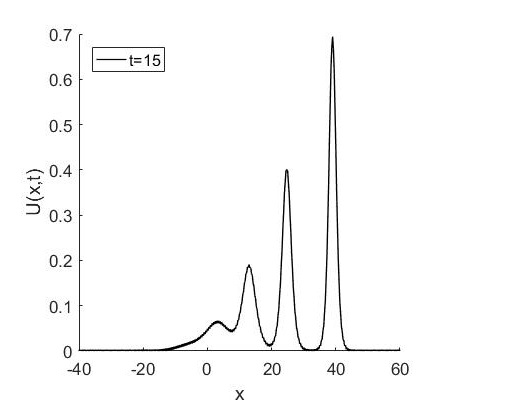}
   \label{fig:3d}
 }
 \caption{Wave generation from an initial positive pulse}
\end{figure}

\noindent
The initial pulse of height $0.4305$ is positioned at $x=5$, Fig \ref{fig:3a}. This pulse is expected to generate new waves immediately at the back as it propagates to the right along the horizontal axis. When the time reaches $t=5$, the height of the frontier is measured as $0.6568$ and it is positioned at $x=18.25$ with an observable distinct follower, Fig \ref{fig:3b}. The follower wave of height $0.3318$ is positioned at $x=12$. Moreover, a second follower, whose bulge is of height $0.1678$ and whose peak is positioned at $x=6.75$, begins to occur at the left of the first follower.
The peak of the frontier reaches at $x=28.75$ as its height is measured as $0.6871$ at $t=10$, Fig \ref{fig:3c}. The shape of the first follower of height $0.3913$ is formed completely and the peak moves to the right along the horizontal axis and positions at $x=18.25$. The second follower of height $0.1736$ positioned $x=9.75$ is also clearly observable with a bulge at the left of it. This bulge is an indicator of a new wave about to come out. When the time is $t=15$, the frontier of height $0.6941$ is positioned at $x=39$ and is separated completely from the first follower, Fig \ref{fig:3d}. The height of the first follower reaches $0.3998$. The peak position of this wave is determined as $x=24.75$. The first follower also separates from the second follower completely at this time. The second follower positioned at $x=13$ reaches $0.1910$ units height. Even though the third follower is not completely formed, the shape of it tends to resemble a solitary wave. One should note that both frontier and the follower waves increase their heights rapidly till they are separated from the previous one. When they are about to be separated completely from the other waves, the increase of their heights keeps but slightly. The velocity of the frontier decreases as new waves are generated as the velocities of both the first and the second followers increase slightly.

\noindent
The conservation laws are determined analytically as
\begin{equation}
{\scriptsize
\begin{aligned}
M_0&=20\,{\it arctanh} \left(  A ^{-1}
 \right) \sqrt {2} \\
E_0&= {\frac {100}{3}}\,{\frac {-A+16\,{\it arctanh}
 \left(  A ^{-1} \right) +4\,{\it 
arctanh} \left(  A^{-1} \right) 
\sqrt {2}\sqrt {7}}{A}}\\
H_0&= -{\frac {4900}{81}}\,{\frac {2114448\,{\it arctanh} \left(  A ^{-1} \right) +565110\,{\it arctanh}
 \left( A ^{-1} \right) \sqrt {2}
\sqrt {7}-148005\,\sqrt {7}-276892\,\sqrt {2}}{13455\,\sqrt {7}+25172
\,\sqrt {2}}}
\end{aligned}
}
\end{equation} 
where $A=(2\sqrt{2}+\sqrt{7})$. These values are computed as $M_0=5.225504574$, $E_0=1.503363838$ and $H_0=1.599480484$ by using symbolic computational mathematical software. The absolute relative changes of the conservation laws at some specific times during the simulation process are reported in Table \ref{t5}. Both polynomial and trigonometric methods seem successful to keep all conservation law values with almost the same decimal digit sensitivity when compared with each other.  
\begin{table}[hp]
\caption{The absolute relative changes of the conservation laws for the wave generation model} 
	\begin{tabular}{llll}
			\hline\hline	
$t$ &  $C(M_{t})$ & $C(E_{t})$ & $C(H_{t})$  \\ 
			\hline\hline	
		& Polynomial Method &&\\
			\hline
$5$  & $8.0719 \times 10^{-7}$ &$3.0588\times 10^{-5}$ &$ 1.2886\times 10^{-3}$  \\ 
$10$  & $2.7652\times 10^{-6}$ &$4.1342\times 10^{-5} $&$ 1.8485\times 10^{-3}$ \\ 
$15$ & $7.0383\times 10^{-6}$ & $6.1132\times 10^{-4}$ & $2.1571\times 10^{-3}$  \\ 
 \hline
& Trigonometric Method &&\\
\hline
$5$  & $2.0371 \times 10^{-7}$ &$ 1.3634\times 10^{-5}$ &$ 9.3847\times 10^{-4}$  \\ 
$10$  & $5.3069\times 10^{-7}$ &$ 1.3773\times 10^{-5} $&$ 1.5776\times 10^{-3}$ \\ 
$15$ & $2.5527\times 10^{-6}$  &$ 5.9309\times 10^{-4}$ &$1.4303\times 10^{-3}$  \\ 
\hline\hline
	\end{tabular}
	\label{t5}
\end{table}

\subsection{Interaction of two Positive Solitary Waves Moving in the Opposite Direction}
\begin{figure}[hp]
    \subfigure[Initial data]{
   \includegraphics[scale =0.6] {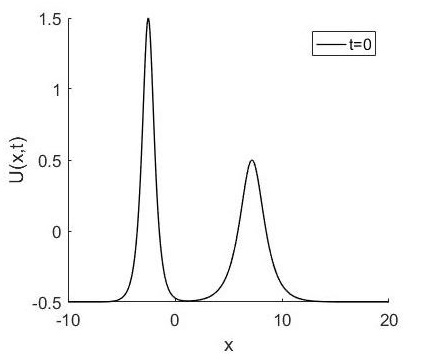}
   \label{fig:4a}
 }
   \subfigure[$t=2$]{
   \includegraphics[scale =0.6] {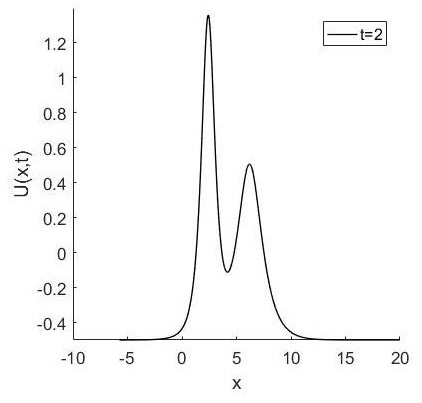}
   \label{fig:4b}
 }
 \subfigure[$t=2.5$]{
   \includegraphics[scale =0.6] {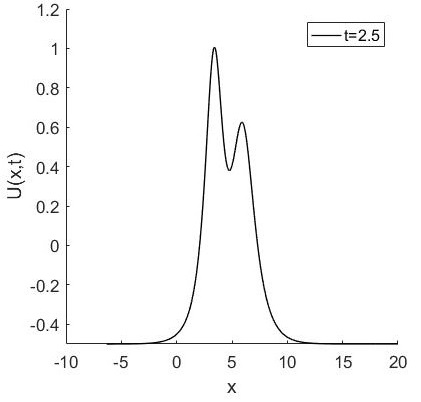}
   \label{fig:4c}
 }
  \subfigure[$t=4$]{
   \includegraphics[scale =0.6] {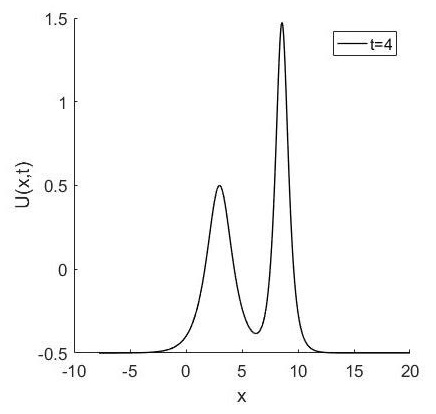}
   \label{fig:4d}
 }
  \subfigure[$t=5$]{
   \includegraphics[scale =0.6] {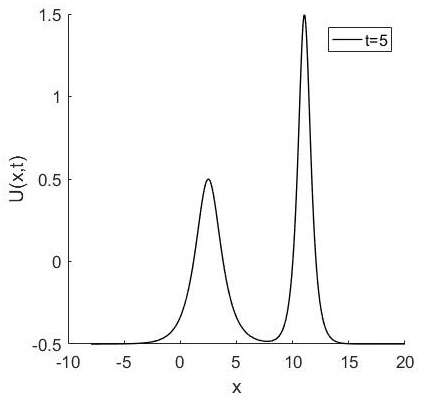}
   \label{fig:4e}
 }
 \caption{Interaction of two positive bell shape solitaries}
\end{figure}

The interaction of two positive bell shape solitaries are also studied by both proposed methods. The initial condition 
\begin{equation}
u(x,0)=-\dfrac{1}{2}+2\dfrac{\left(e^{x-5}+2e^{2x+5}\right)\left(1-\dfrac{1}{9}e^{3x}\right)+\dfrac{1}{3}e^{3x}\left(e^{x-5}+2e^{2x+5}\right)}{\left(e^{x-5}+e^{2x+5}\right)^2+\left(1-\dfrac{1}{9}e^{3x} \right)^2}
\end{equation}
is a particular form the $N$-soliton solution derived from the solution given in \cite{wazwazbook}. This initial condition gives two well seperated positive bell shaped solitaries of heights $1.499631216748184$ and $0.499996193761888$ positioned at $x=-2.5$ and $x=7.2$, respectively, at the beginning, Fig \ref{fig:4a}. Both solitaries propagate in the ooposite directions along the horizontal axis as time goes.

\noindent
We choose the parameters $\alpha=6$, $\beta=6$ and $\mu=1$ in the Gardner equation. The designed routines are run up to the terminating time $t=5$ with the discretization parameters $N=600$ and $\Delta t=0.01$ in the finite problem interval $[-10,20]$.

\noindent
When the time reaches $t=2$, it is observed that the interaction has started, Fig \ref{fig:4b}. The height of the higher solitary is measured as $1.3598$ and its peak is positioned as $x=2.4$. The height of the lower one increases to $0.5067$ and the position of its peak is determined at $x=6.2$. The height of the higher wave reaches $1.0057$ as the height of the lower one $0.6252$ at the time $t=2.5$, Fig \ref{fig:4c}. The peaks of both the higher and the lower solitaries are positioned at $x=3.4$ and $x=5.90$, respectively. When the time reaches $t=4$, the solitaries begin to separate, Fig \ref{fig:4d}. The height of the higher one increases to $1.4732$ and it is positioned at $x=8.550$. The peak of the lower one is positioned at $x=3$ and the height of it decreased to $0.4999$. At the end of the simulation, we observe both solitaries are separated completely and return to their original shapes and heights, Fig \ref{fig:4e}. The heights of both solitaries are determined as $0.4997$ and $1.4959$ as the peaks reach $x=2.5$ and $11.05$ as keeping to propagate on their own ways.

\noindent
The conservation laws for the interaction of the solitary waves are also computed during the simulations for both routines. The values are computed as $M_0=-8.716821423$, $E_0=7.216821423$ and $H_0=-2.34182152$ initially. The absolute relative changes of the conservation laws are tabularised in Table \ref{t6}.  
\begin{table}[hp]
\caption{The absolute relative changes of the conservation laws for the interaction of two positive bell shaped solitaries propagating in the opposite direction} 
	\begin{tabular}{llll}
			\hline\hline	
$t$ &  $C(M_{t})$ & $C(E_{t})$ & $C(H_{t})$  \\ 
			\hline\hline	
		& Polynomial Method &&\\
			\hline
$1$  & $1.0880 \times 10^{-5}$ &$4.6286\times 10^{-4}$ &$ 2.8376\times 10^{-4}$  \\
$2.5$  & $8.5942\times 10^{-6}$ &$1.6490\times 10^{-4} $&$ 3.5258\times 10^{-3}$ \\

$5$ & $2.5123\times 10^{-4}$ & $4.0562\times 10^{-4}$ & $1.4966\times 10^{-4}$  \\ 
 \hline
& Trigonometric Method &&\\
\hline
$1$  & $1.1763 \times 10^{-5}$ &$ 4.7382\times 10^{-5}$ &$ 2.5578\times 10^{-4}$  \\ 
$2.5$  & $8.6367\times 10^{-6}$ &$ 1.6517\times 10^{-4} $&$ 2.6099\times 10^{-3}$ \\

$5$ & $1.9378\times 10^{-4}$  &$ 3.3646\times 10^{-4}$ &$6.0904\times 10^{-5}$  \\ 
\hline\hline
	\end{tabular}
	\label{t6}
\end{table}

\section{Conclusion}
\noindent
The collocation methods based on cubic polynomial and trigonometric B-spline functions are derived for the numerical solutions of some analytical and non-analytical problems for the Gardner equation. Having no continuous derivatives of both sets of basis functions force us to reduce the order of the dispersion term $u_{xxx}$ to lead a coupled system of equations. Assuming the solution and its derivative in the same space satisfies to approximate to both functions by using the same basis function set. The time integration of the system is achieved by the classical Crank-Nicolson method owing to its large stability region. We measure the errors between the numerical and the analytical solutions in case the existence of the analytical solutions. The lowest three conservation laws representing momentum, energy and the Hamiltonian quantity are computed and the relative absolute changes of all laws are examined in details even for non-analytical solutions for a fair comparison of both methods.

\noindent
In the first example, we investigate the numerical solution describing the propagation of a bell shaped positive solitary along the horizontal axis. Various discretization parameters are used to obtain the numerical solutions. Even though the results obtained by the polynomial method seems better when the number of grids are less, the increase of the grid numbers improved the solutions more in the trigonometric method. The comparison of absolute relative changes of conservation laws does not give a clear idea about the better method since the results obtained by both methods pretty accurate.

\noindent
In the kink type solution of the Gardner equation, the polynomial method gives more accurate results for all choices of number of grids. The conservation laws obtained by both methods are sensitive in the same decimal digits in almost all cases.

\noindent
The perturbation of a single positive bell shaped solitary wave is derived to study wave generation for the Gardner equation successfully. Both methods generate valid solutions. The sensitivities of the absolute relative changes of the lowest conservation laws are indicators of two efficient methods.

\noindent
Finally, we study the interaction of two positive bell shaped solitaries propagating in the opposite directions. Both methods simulate the solutions successfully. The conservation laws are also kept during the full elastic collision and later.
 
\noindent
\textbf{\textit{Acknowledgements:}} \textsl{This study is a part of the project with number 2016/19052 supported by Eskisehir Osmangazi University Scientific Research Projects Committe, Turkey.}

\end{document}